\documentclass[12pt]{amsart}

\usepackage{amssymb}
\usepackage{latexsym}
\usepackage{verbatim}
\usepackage{fullpage}
\usepackage{pslatex}

\input {cyracc.def}
 \font\tencyr=wncyr10 
\font\tencyi=wncyi10 
\font\tencysc=wncysc10 
\def\rus{\tencyr\cyracc}
\def\rusi{\tencyi\cyracc}
\def\rusc{\tencysc\cyracc}

\renewenvironment{proof}
{\noindent {\sl Proof.}\quad }{\hfill
$\square$ \vskip1.1ex\noindent }

\newenvironment{proof*}
{\noindent {\sl Proof.}\quad }{\hfill
$\square$}

\renewcommand{\theequation}{\thesection .\arabic{equation}}
\renewcommand{\thesubsubsection}{\theequation .\arabic{subsubsection}}

\catcode`\@=\active
\catcode`\@=11
\def\@eqnnum{\hbox to
.01pt{}\rlap{\hskip-\displaywidth(\mathbf{\theequation})}}
\catcode`\@=12

\newenvironment{s}[1]
{ \vskip1.2ex \refstepcounter{equation}
\noindent {\bf \theequation\quad #1.} \begin{sl}}{\end{sl}
\vskip1.1ex\noindent }

\newenvironment{rem}[1]
{ \vskip1.2ex \refstepcounter{equation}
\noindent {\bf \theequation\enspace {#1}.} }{ \vskip1.1ex\noindent }

\newenvironment{subs}[1]
{\vskip1.2ex \refstepcounter{equation}
\noindent {\bf (\theequation)\quad #1.} }{\quad}

\newcommand {\ah}{{\frak a}}
\newcommand {\be}{{\frak b}}

\newcommand {\g}{{\frak g}}
\newcommand {\h}{{\frak h}}

\newcommand {\te}{{\frak t}}

\newcommand {\esi}{\varepsilon}
\newcommand {\ap}{\alpha}

\newcommand {\lb}{\lambda}

\newcommand {\vp}{\varphi}

\newcommand {\ca}{{\mathcal A}}

\newcommand {\MM}{{\mathcal M}}

\newcommand {\VV}{{\Bbb V}}
\newcommand {\W}{{\Bbb W}}

\newcommand {\rk}{{\mathrm{rk\,}}}

\newcommand {\GR}[2]{{\textrm{{\bf #1}}}_{#2}}

\newcommand {\Ab}{{\frak Ab}}

\newcommand {\beq}{\begin{equation}}
\newcommand {\eeq}{\end{equation}}

\newcommand{\vts}{\VV(\theta_s)}

\newfam\Bbbfam\newfam\eufam%
\font\Bbbfont=msbm10 scaled 1200%
\font\olala=msam10 scaled 1200%
\font\frak=eufm10 scaled 1400%
\font\Bbbsmallfont=msbm8%
\textfont\Bbbfam=\Bbbfont\scriptfont\Bbbfam=\Bbbsmallfont
\font\euzw=eufm10 scaled 1200%
\font\euac=eufm7 scaled 1200%
\font\euacc=eufm7 scaled 1000%
\textfont\eufam=\euzw\scriptfont\eufam=\euac%
\scriptscriptfont\eufam=\euacc%
\def\frak{\fam\eufam}%
\def\Bbb{\fam\Bbbfam}%

\def\varnothing{\hbox {\Bbbfont\char'077}}
\def\square{\hbox {\olala\char"03}}

\begin{document}
\setlength{\parskip}{2pt plus 4pt minus 0pt}
\hfill {\scriptsize March 18, 2003} 
\vskip1ex
\vskip1ex

\title[]{Long Abelian ideals
}
\author[]{\sc Dmitri I. Panyushev}
\thanks{This research was supported in part by RFBI Grants no.
01--01--00756 and  02--01--01041}
\maketitle
\begin{center}
{\footnotesize
{\it Independent University of Moscow,
Bol'shoi Vlasevskii per. 11 \\
121002 Moscow, \quad Russia \\ e-mail}: {\tt panyush@mccme.ru }\\
}
\end{center}

Let $\be$ be a Borel subalgebra of a simple Lie algebra $\g$.
Let ${\Ab}$ denote the set of all Abelian ideals of $\be$.
It is easily seen that any $\ah\in \Ab$ is actually contained in the nilpotent
radical of $\be$. Therefore $\ah$ is determined by the the corresponding
set of roots. More precisely, let $\te$ be a Cartan subalgebra of $\g$
lying in $\be$ and let $\Delta$ be the root system of the pair $(\g,\te)$.
Choose the system of positive roots,
$\Delta^+$, so that the roots of $\be$ are positive.
Then $\ah=\oplus_{\gamma\in I}\g_\gamma$, where $I$ is a suitable subset
of $\Delta^+$
and $\g_\gamma$ is the root space for $\gamma\in\Delta^+$.

A nice result of D.\,Peterson says that the cardinality of
$\Ab$ is $2^{\rk\g}$. His approach uses a one-to-one correspondence between
the Abelian ideals and the so-called `minuscule' elements of the affine
Weyl group $\widehat W$ (see Section~\ref{prelim} for precise definitions).
 An exposition of Peterson's results is found in \cite{Ko1}. Peterson's work
appeared to be the point of departure for active recent investigations of
Abelian ideals, 
and related problems of
representation theory and combinatorics 
\cite{CP},\cite{CP2},%
\cite{Ko1},%
\cite{lp},\cite{pr}. Our definition of minuscule elements follows
Kostant's paper \cite{Ko1}, so that $w\in\widehat W$ is minuscule in
our sense if and only if $w^{-1}$ is minuscule in the sense of Cellini--Papi
\cite{CP},\,\cite{CP2}.
An elegant proof of Peterson's theorem is given in \cite{CP}.
Let $\ca$ be the fundamental alcove of $\widehat W$.
Cellini and Papi show that $w\in \widehat W$ is minuscule
if and only if $w^{-1}{\cdot}\ca\subset 2\ca$. Since $2\ca$ consists of
$2^{\rk\g}$ alcoves and $\widehat W$ acts simply transitively on the
set of alcoves, Peterson's theorem follows.

In this paper, we first show that methods of \cite{CP} can be adapted to
solving the following problems:

{\it Suppose\/ $\g$ has two root lengths.
\begin{itemize}
\item describe (enumerate) the Abelian ideals such that the
corresponding set\/ $I$ consists only of long roots (such Abelian ideals
are said to be {long});
\item give a characterization of the corresponding (= long)
minuscule elements.
\end{itemize}
}

\noindent
Write $\Ab_l$ (resp. $\MM_l$) for the set of long Abelian ideal
(resp. long minuscule elements).
We give a uniform answer to both problems. Let $d>1$ be the ratio of squares
of root lengths and $a$ the number of long simple roots. Then
\[
   \#(\Ab_l)=d^a \ .
\]
To obtain a characterisation of the long minuscule elements, we consider
a certain simplex $\ca_s$ that lies between $\ca$ and  $2\ca$.
Then the answer is that $w\in \MM_l$ if and only if $w^{-1}{\cdot}\ca\subset
\ca_s$ if and only if any reduced decomposition of $w$ contains no short
simple reflections (we say that a reflection is short, if it corresponds
to a short root). We also describe the rootlets of long Abelian ideals.
\\
Second, we show that the theory of long Abelian ideals is closely related to
describing of commutative subalgebras in the isotropy representations of
some symmetric spaces. Let $\theta$ be the highest root
and $\theta_s$ the short dominant root of $\Delta^+$. Write $\VV(\lb)$
for the $\g$-module with highest weight $\lb$.
Assume that $|\theta|^2/|\theta_s|^2=2$, i.e., we exclude the case of
$\GR{G}{2}$. Then $\tilde\g=\g\oplus\vts$ has a structure of
a ${\Bbb Z}_2$-graded Lie algebra, so that it makes sense to speak
about commutative subalgebras of $\vts$.
In Section~\ref{vee}, we prove that there is a natural one-to-one
correspondence between the $\be$-stable commutative subalgebras
of $\vts$ and the long Abelian $\be^\vee$-ideals in
$\g^\vee$. Here $\g^\vee$ is the Langlands dual Lie algebra for $\g$.
\\[.6ex]
{\bf Acknowledgements.} This paper was written during my stay at the
Max-Planck-Institut f\"ur Mathematik (Bonn). I would like to thank this
institution for the warm hospitality and support.

\section{Notation and other preliminaries} \label{prelim}

\noindent
\begin{subs}{Main notation}
\end{subs}
$\Delta$ is the root system of $(\g,\te)$ and
$W$ is the usual Weyl group. For $\ap\in\Delta$, $\g_\ap$ is the
corresponding root space in $\g$.

$\Delta^+$ \ is the set of positive
roots and $\rho=\frac{1}{2}\sum_{\ap\in\Delta^+}\ap$.

$\Pi=\{\ap_1,\dots,\ap_p\}$ \ is the set of simple roots in $\Delta^+$.

$\vp_1,\dots,\vp_p$  \ are the fundamental weights corresponding to $\Pi$.
 \\
 We set $V:=\te_{\Bbb R}=\oplus_{i=1}^p{\Bbb R}\ap_i$ and denote by
$(\ ,\ )$ a $W$-invariant inner product on $V$.   As usual,  $\mu^\vee=2\mu/(\mu,\mu)$ is the coroot
for $\mu\in \Delta$.
Letting $\widehat V=V\oplus {\Bbb R}\delta\oplus {\Bbb R}\lb$, we extend
the inner product $(\ ,\ )$ on $\widehat V$ so that $(\delta,V)=(\lb,V)=(\delta,\delta)=
(\lb,\lb)=0$ and $(\delta,\lb)=1$.

$\widehat\Delta=\{\Delta+k\delta \mid k\in {\Bbb Z}\}$ is the set of
affine real roots
and $\widehat W$ is the  affine Weyl group.
\\
Then $\widehat\Delta^+= \Delta^+ \cup \{ \Delta +k\delta \mid k\ge 1\}$ is
the set of positive
affine roots and $\widehat \Pi=\Pi\cup\{\ap_0\}$ is the corresponding set
of affine simple roots.
Here $\ap_0=\delta-\theta$, where $\theta$ is the highest root
in $\Delta^+$.  The inner product $(\ ,\ )$ on $\widehat V$ is
$\widehat W$-invariant.
\\
For $\ap_i$ ($0\le i\le p$), we let $s_i$ denote the corresponding simple reflection in $\widehat W$.
The length function on $\widehat W$ with respect
to  $s_0,s_1,\dots,s_p$ is denoted by $\ell$.
We consider two different actions of $\widehat W$: the linear action on
$\widehat V$ and the affine-linear action on $V$. To distinguish these two, we use
dot `$\cdot$' for denoting the second action.
For any $w\in\widehat W$, we set
\[
\widehat N(w)=\{\ap\in\widehat\Delta^+ \mid w(\ap) \in -\widehat \Delta^+ \} .
\]
Our convention concerning $\widehat N(w)$ is the same as in
\cite{Ko1},\,\cite{lp}
but is opposite to that in \cite{CP},\,\cite{CP2}.
\begin{subs}{Abelian ideals}
\end{subs}
Let $\ah\subset\be$ be an Abelian ideal. It is easily seen that
$\ah\subset [\be,\be]$. Therefore $\ah=\underset{\ap\in I}{\oplus}\g_\ap$
for a subset $I\subset \Delta^+$, which is called the {\it set of roots of\/}
$\ah$.
In what follows, an Abelian ideal
will be identified with the respective set of roots.
That is, $I$ is said to be an Abelian ideal, too.  We
also say that $\ah$ is a {\it geometric\/}
Abelian ideal, while $I$ is a {\it combinatorial\/} Abelian ideal.
In the combinatorial context, the definition of an Abelian ideal
can be stated as follows.
\\
$I\subset\Delta^+$ is an Abelian ideal, if the following two conditions
are satisfied:

(a) for any $\mu,\nu\in I$, we have $\mu+\nu\not\in\Delta$;

(b) if $\gamma\in I$, $\nu\in\Delta^+$, and $\gamma+\nu\in\Delta$, then
$\gamma+\nu\in I$.
\\[.6ex]
Following D.\,Peterson, an element $w\in \widehat W$ is said to be
{\it minuscule},
if $\widehat N(w)$ is of the form $\{\delta-\gamma \mid \gamma\in I \}$,
where $I$ is a subset of $\Delta^+$.
It was shown by Peterson
that such an $I$ is a combinatorial Abelian
ideal and, conversely, each Abelian ideal occurs in this way,
see \cite[Prop.\,2.8]{CP},\,\cite{Ko1} .
Hence one obtains a one-to-one correspondence between the
Abelian ideals of $\be$ and the minuscule elements of $\widehat W$.
The set of minuscule elements in $\widehat W$
is denoted by $\MM$.
If $w\in\MM$, then we write $I_w$ (resp. $\ah_w$) for the corresponding
combinatorial
(resp. geometric) Abelain ideal. That is,
\[
   I_w=\{\gamma\in\Delta^+ \mid \delta-\gamma \in\widehat N(w)\}
\ \textrm{ and }\  \ah_w=\oplus_{\ap\in I_w}\g_\ap \ .
\]
Conversely, given $I\in \Ab$, we write $w\langle I\rangle$ for the
respective minuscule element. Notice that
\[
  \dim\ah_w =\# (I_w)=l(w) \ .
\]
Throughout the paper, $I$ or $I_w$ stands for a combinatorial Abelian ideal.

\section{Long minuscule elements and long Abelian ideals}
\setcounter{equation}{0}

\noindent
From now on, we assume that  $\Delta^+$ has two root lengths.
To distinguish different objects related to
long and short roots, we use the subscripts `$l$' and `$s$'.
For instance,
$\Pi_l$ is the set of long simple roots
and $\Delta^+_l$ is the set of long positive roots. Accordingly, each
simple reflection $s_i$ is either short or long. Since $\theta$ is long,
the reflection $s_0$ is regarded as long.
Write $\theta_s$ for the unique short dominant root in $\Delta^+$.

We say that $I\in\Ab$ is \textit{long}, if
$I\subset \Delta^+_l$.
Write $\Ab_l$ for the set of all long Abelian ideals. The corresponding
subset of $\MM$ is denoted by $\MM_l$.
Notice that the analogous notion of a short Abelian ideal does not make sense,
for any non-empty Abelian ideal contains the highest root $\theta$, which is
long.

\begin{s}{Proposition}   \label{charact1}
Suppose $w\in\MM$.
Then
$w\in \MM_l$ if and only if any reduced decomposition of $w$
does not contain short simple reflections.
\end{s}\begin{proof}
`$\Rightarrow$' \ Suppose $w=w''s_iw'$, where $\ap_i$ is a short simple root.
As was noticed in \cite[Sect.\,2]{lp}, every right substring in
a reduced decomposition of $w$ is again a minuscule
element and the corresponding (combinatorial) ideal is a subset of $I_w$.
In particular, we have $s_iw', w'\in\MM$. Furthermore,
\[
   I_{s_iw'}=\{\gamma\}\cup I_{w'} ,
\]
where $w'(\delta-\gamma)=\ap_i$. It follows that $\gamma$ is a short root
lying in $I_{s_iw'}\subset I_w$.

`$\Leftarrow$' \ Agrue by induction on the length of $w$.
If $w'\in \MM_l$, $w\in \MM$, and $w=s_iw'$ for a long reflection $s_i$, then
$I_w=\{\gamma\}\cup I_{w'}$, where $w'(\delta-\gamma)=\ap_i$.
Since $\gamma$ is long, we conclude that $w\in\MM_l$.
\end{proof}%
Consider the collection of affine hyperplanes in $V$
\[
    H_{\mu,k}=\{x\in V \mid (\mu, x)=k\} \ ,
\]
where $\mu\in\Delta^+$ and $k\in {\Bbb Z}$.
The connected components of
$\displaystyle V \setminus \bigcup_{\mu,k}
H_{\mu,k}$ are called {\it alcoves\/}.
It is well-known that $\widehat W$, as group of affine transformations of
$V$, acts simply-transitively on the set of
alcoves.
Recall that the fundamental alcove of $\widehat W$ is
\[
   \ca=\{x\in V \mid (\ap,x)>0 \ \forall \ap\in\Pi \ \& \ (\ap,\theta)<1\} \ .
\]
Set
\[
   \ca_s=\{x\in V \mid (\ap,x)>0 \ \forall \ap\in\Pi \ \& \ (\ap,\theta_s)<1\}
  \ .
\]
It is clear that $\ca_s$ is a union of (finitely many) alcoves and it contains
$\ca$.

\begin{s}{Proposition}   \label{charact2}
Let $w\in\widehat W$ be a minuscule element.
We have
$I_w\in \Ab_l$ if and only if $w^{-1}{\cdot}\ca\subset \ca_s$.
\end{s}\begin{proof}
We use the following result proved in \cite[1.1]{CP}:
\begin{align*}
   -\mu+k\delta\in N(w) \ & \Leftrightarrow \ H_{\mu,k} \text{
\ \ separates $\ca$ and $w^{-1}{\cdot}\ca$}\qquad (k > 0) \\
 \mu+k\delta\in N(w) \ &\Leftrightarrow \ H_{\mu,-k} \text{
separates $\ca$ and $w^{-1}{\cdot}\ca$} \qquad (k\ge 0) \ .
\end{align*}
We also know that $w^{-1}{\cdot}\ca\subset 2\ca$, since $w\in\MM$.

(a) If $w^{-1}{\cdot}\ca\not\subset \ca_s$, then the hyperplane
$H_{\theta_s,1}$ separates $\ca$ and $w^{-1}{\cdot}\ca$. Hence
$\delta-\theta_s\in N(w)$, i.e., $\theta_s\in I_w$.

(b) Conversely, if $w^{-1}{\cdot}\ca\subset \ca_s$, then
any hyperplane $H_{\mu,k}$ separating $\ca$ and $w^{-1}{\cdot}\ca$ must meet
$\ca_s$. If $x\in \ca_s\cap H_{\mu,k}$ and $\mu\in \Delta^+_s$, then
$0< (x,\mu)\le (x,\theta_s) < 1$, which is impossible.
Hence the hyperplanes $H_{\mu,1}$, with $\mu$ short, do not separate
$\ca$ and $w^{-1}{\cdot}\ca$. Thus, $I_w\cap \Delta^+_s=\varnothing$.
\end{proof}%
It follows that the number of long Abelian ideals is equal to the
number of alcoves that fit in $\ca_s$. In other words,
\begin{equation}  \label{abl}
  \#(\Ab_l) = \mathrm{vol\,}(\ca_s)/ \mathrm{vol\,}(\ca) \ .
\end{equation}
\vskip-1.2ex
\begin{s}{Theorem}  \label{volume}  \par
\hbox to \textwidth{
\hfil $
  \displaystyle
\frac{\mathrm{vol\,}(\ca_s)}{\mathrm{vol\,}(\ca)} =
  \left(\frac{|\theta|^2}{|\theta_s|^2}\right)^{\#\Pi_l}$.
\hfil }
\end{s}\begin{proof*}
Write $\theta=\sum_{i=1}^p m_i \ap_i$ and
$\theta_s=\sum_{i=1}^p c_i \ap_i$. Then
\[
 \frac{\mathrm{vol\,}(\ca_s)}{\mathrm{vol\,}(\ca)} =
  \prod_{i=1}^p \frac{(\vp_i,\theta)}{(\vp_i,\theta_s)}
 =\prod_{i=1}^p \frac{m_i}{c_i} \ .
\]
Consider the dual root system $\Delta^\vee$. We have
$\{\ap^\vee \mid \ap\in\Pi\}$ is a set of simple roots,
$\theta^\vee_s$ is the highest root in $\Delta^\vee$, and
\[
 \theta^\vee_s=\frac{2}{(\theta_s,\theta_s)}\sum_{i=1}^p c_i \ap_i=
\sum_{\ap_i\in \Pi_s}c_i\ap_i^\vee + \sum_{\ap_i\in \Pi_l}c_i
\frac{|\theta|^2}{|\theta_s|^2}\ap_i^\vee \ .
\]
On the other hand, the collection of the coefficients of the
highest root in $\Delta^\vee$ is the same as in $\Delta$.
Hence
\[
  \prod_{i=1}^p m_i=\prod_{i=1}^p c_i{\cdot}
\left(\frac{|\theta|^2}{|\theta_s|^2}\right)^{\#\Pi_l} \ ,
\]
and we are done.
\end{proof*}%
\begin{s}{Corollary}
The number of long Abelian ideals equals:

$2^{p-1}$ for ${\frak so}_{2p+1}$; \
2 \ for ${\frak sp}_{2p}$; \
4  \ for $\GR{F}{4}$; \
3  \ for $\GR{G}{2}$.
\end{s}%
In \cite{lp}, we introduced the notion of the {\it rootlet\/} of an
Abelian ideal.  Let $I_w$, $w\in\MM$, be a non-trivial Abelian $\be$-ideal.
Then $\tau(I_w):=w(\ap_0)+\delta$ is a long positive root \cite[2.4]{lp},
which is called the {\it rootlet\/} \ of $I_w$. It is natural to inquire as
to what the rootlets of long Abelian ideals are.

\begin{s}{Proposition}   \label{rootlet}
Suppose $I$ is a nontrivial Abelian ideal. Then
$I\in \Ab_l$ if and only if $\theta - \tau(I)$ is a
certain sum of only long simple
roots, i.e., $\theta$ and $\tau(I)$ have equal coefficients of
all short simple roots.
\end{s}\begin{proof}
`$\Rightarrow$' \quad
By Proposition~\ref{charact1}, a reduced decomposition of
$w=w\langle I\rangle$
contains only reflections $s_0$ and $s_i$, $i>0$, for $\ap_i\in\Pi_l$.
We argue by induction on $\ell(w)$.

(a) If $\ell(w)=1$, then
$w=s_0$ and $s_0(\ap_0)+\delta=\theta$.

(b) Suppose $w\in\MM_l$ and $w(\ap_0)+\delta=\gamma\in\Delta^+_l$.
If $s_iw\in \MM$, then $s_i(\gamma-\delta)+\delta=s_i(\gamma)$.
Consider two possibilities for $s_i$.
If $s_i$ corresponds to $\ap_i\in\Pi_l$, then $\gamma$ and $s_i(\gamma)$
have equal coefficients of the short simple roots.
If $s_i=s_0$, then
\[
   s_0(\gamma)=\left\{ \begin{array}{cl}
    \gamma, & \text{ if $(\gamma,\theta)=0$} \\
   \delta-(\theta-\gamma), &  \text{ if $(\gamma,\theta)\ne 0$ .}
    \end{array}\right.
\]
As the last case is impossible in the minuscule situation (for,
we would obtain a root which is not in $\Delta^+$),
the induction step is complete.

`$\Leftarrow$' \quad  Assume that a reduced decomposition of $w$
contains a short reflection. Then the first occurrence of
short reflections will certainly reduce some short simple root
coefficients of the current rootlet.
It is conceivable that a consequent occurrence
of short reflections would restore the previous coefficients of the short
simple roots. However, it is not possible in view of
\cite[Corollary\,3.3]{lp}. Thus, the reduced decompositions of $w$ cannot
contain short simple reflections.
\end{proof}%
In \cite{lp}, we studied the poset structure of the fibres of the mapping
$\tau: \Ab\setminus\{\varnothing\} \to \Delta^+_l$ that takes a non-trivial
Abelian ideal to its rootlet. By the previous proposition, we have the
set of non-trivial long Abelian ideals is the union of fibres of this mapping.

\begin{rem}{Examples} We use the notation and numbering of roots as in
\cite[Tables]{VO}.\\
1. $\g={\frak sp}_{2p}$.
Here $\theta=2\esi_1=2\vp_1$ and $\theta_s=\esi_1+\esi_2=\vp_2$.
Normalize the bilinear form $(\ ,\ )$ so that
$(\esi_i,\esi_i)=1$. Then
$\ca$ has the vertices $\esi_1,\esi_1+\esi_2,\dots,
\esi_1+\esi_2+\ldots +\esi_p$, whereras $\ca_s$ has the vertices
$2\esi_1,\esi_1+\esi_2,\dots,
\esi_1+\esi_2+\ldots +\esi_p$. The unique non-trivial long
Abelian ideal is $\{\theta\}$.

2. $\g={\frak so}_{2p+1}$. Here $\theta=\esi_1+\esi_2=\vp_2$ and
$\theta_s=\esi_1=\vp_1$.
The unique maximal long Abelian ideal consists of
the roots $\{\esi_i+\esi_j \mid 1\le i< j\le p\}$. The only generator is
$\esi_{p-1}+\esi_p$.
The corresponding rootlet is also $\esi_{p-1}+\esi_p=\ap_{p-1}+2\ap_p$.

3. $\g=\GR{F}{4}$.  Here $\theta=2\ap_1+4\ap_2+3\ap_3+2\ap_4=\vp_4$ and
$\theta_s=2\ap_1+3\ap_2+2\ap_3+\ap_4=\vp_1$.
The non-trivial long Abelian ideals appear in the first three
rows of Table~1 in \cite{lp}. The unique maximal long Abelian ideal is
$\{\theta,\theta-\ap_4,\theta-\ap_4-\ap_3\}$.

4. $\g=\GR{G}{2}$. The non-trivial long Abelian ideals are
$\{3\ap_1+2\ap_2\},\,\{3\ap_1+2\ap_2,3\ap_1+\ap_2\}$.
\end{rem}

\section{Long abelian ideals and little adjoint representations}
\setcounter{equation}{0}   \label{vee}

\noindent
In this section, we elaborate on a relationship between long Abelian ideals
and commutative subalgebras in the isotropy representations
of some symmetric spaces.

We still assume that $\g$ has two root lengths. Then the $\g$-representation
with highest weight $\theta_s$ is said to be {\it little adjoint\/}.
The corresponding $\g$-module is denoted by $\vts$.
The properties of the $\g$-module $\vts$ are similar to that of $\g=\VV(\theta)$.
For instance, the set of nonzero weights of $\vts$ is
$\Delta_s$; the dimension of the zero-weight space is
$\#(\Pi_s)$ and each weight space corresponding to $\mu\in\Delta^+_s$ is
one-dimensional, see \cite[2.8]{ya-tg}.
\\[.6ex]
From now on, we stick to the case, where
$\frac{|\theta|^2}{|\theta_s|^2}=2$. Then $\vts$ is the isotropy representation of
a symmetric space; more precisely,
\[
\tilde\g=\g\oplus \vts
\]
has a natural structure of ${\Bbb Z}_2$-graded
simple Lie algebra. The following table shows all possibilities for
$\theta_s$ and $\tilde\g$.
\begin{center}
\begin{tabular}{c||ccc}
$\g$   &  ${\frak so}_{2p+1}$ & ${\frak sp}_{2p}$ & $\GR{F}{4}$ \\ \hline
$\theta_s$ & $\vp_1$ & $\vp_2$ & $\vp_1$ \\ \hline
$\dim\vts$ & $2p+1$ & $2p^2-p-1$ & $26$ \\ \hline
$\tilde\g$ & ${\frak so}_{2p+2}$ & ${\frak sl}_{2p}$ & $\GR{E}{6}$ \\
\end{tabular}
\end{center}
Hence, given a subspace $\ah\subset \vts$, it makes sense to say that it can
be a commutative subalgebra (of $\tilde\g$).

The problem of describing Abelian $\be$-ideals can be generalized to
the setting of ${\Bbb Z}_2$-graded Lie algebras as follows:
\\[.7ex]
{\it Suppose\/ $\hat\h=\h\oplus \VV$ is a
semisimple ${\Bbb Z}_2$-graded Lie algebra,
i.e., $\VV$ is an $\h$-module and $[\VV,\VV]\subset \h$. Let $\be(\h)$ be
a Borel subalgebra of $\h$.
Describe (enumerate) the $\be(\h)$-stable commutative subalgebras in $\VV$.}
\\[.5ex]
The interest of commutative subalgebras in $\VV$ is explained by the fact
that if $\dim\ah=k$, then the corresponding decomposable $k$-vector
in $\wedge^k\VV$ is an eigenvector a Casimir element of $\h$ whose
eigenvalue is maximal possible (and equal $k/2)$; the converse is also true,
see \cite[Sect.\,4]{cas}.
I do not think that an answer to the above problem can be given in
a uniform way (see also examples below). But, for ${\Bbb Z}_2$-gradings
connected with little adjoint representations, there is a reasonably
nice description.

Let $\g^\vee$ denote the Langlands dual Lie algebra for $\g$,
i.e., the root system of $\g^\vee$ is $\Delta^\vee$.
We have $(\Delta^+)^\vee$ is a set of positive roots in $\Delta^\vee$.
Write $\be^\vee$ for the respective Borel subalgebra of $\g^\vee$.
Notice that $(\Delta^+_l)^\vee=(\Delta^\vee)^+_s$ and
$\VV(\theta^\vee)$ is the little adjoint module for $\g^\vee$.
Given a set $S$ of short roots in $\Delta^\vee$, we write
$\ah(S)$ for the respective subspace of $\VV(\theta^\vee)$.

\begin{s}{Lemma}
If $\ah$ is a commutative $\be^\vee$-stable subalgebra of $\VV(\theta^\vee)$,
then $\ah=\ah(S)$ for some $S\subset (\Delta^+_l)^\vee$.
\end{s}\begin{proof}
It follows from \cite[Prop.\,4.9]{cas} that $\ah$ has no zero weight and
$\ah=\ah(S)$ for some $S$ lying in an open halfspace of $V$. However, if
$-\mu^\vee\in S$ for some $\mu\in (\Delta^+_l)^\vee$, then the $\be^\vee$-invariance
of $\ah$ implies that it has a nonzero component in the zero-weight space of
$\VV(\theta^\vee)$.
\end{proof}%
The following result asserts that there is a bijection between
the long Abelian $\be$-ideals in $\g$ and the $\be^\vee$-stable commutative
subalgebras of $\VV(\theta^\vee)$.

\begin{s}{Theorem}   \label{dual} \\
$I\subset\Delta^+_l$ is an Abelian ideal \ $\Leftrightarrow$ \
$\ah(I^\vee)$ is a $\be^\vee$-stable commutative subalgebra in
$\VV(\theta^\vee)$.
\end{s}\begin{proof}
`$\Leftarrow$' \\
1. Suppose $\mu\in I$, $\beta\in \Delta^+$, and
$\mu^\vee+\beta^\vee\in \Delta^\vee$.

-- \ If $\beta$ is long, then $\mu^\vee+\beta^\vee=
(\mu+\beta)^\vee\in I^\vee$.

-- \ If $\beta$ is short, then $\mu+\beta\in I$ is a short root, which is
impossible.\\
Thus, the space $\ah(I^\vee)$ is $\be^\vee$-stable.
\\
2. Let $\mu,\beta\in I$ and assume that
$\mu^\vee+\beta^\vee\in \Delta^\vee$.
Since $I$ is Abelian, we have $\mu+\beta\not\in\Delta$. Then
$\mu+\beta=2\gamma$, where $\gamma\in\Delta_s$. Therefore $(\mu,\beta)=0$
and $(\mu,\gamma)=(\beta,\gamma)>0$. Hence
$\gamma-\mu,\gamma-\beta\in \Delta_s$ and $(\gamma-\mu)+(\gamma-\beta)=0$.
Because one of these roots is positive, we conclude that
$\gamma\in I\cap \Delta_s$. This contradiction proves that
$\mu^\vee+\beta^\vee\not\in \Delta^\vee$, i.e.,
$\ah(I^\vee)$ is a commutative subalgebra.

`$\Rightarrow$' \\
The argument is similar. One has only use the following property of this
${\Bbb Z}_2$-grading:

Suppose $\nu,\gamma\in (\Delta^\vee)_s$ and  $v_\nu, v_\gamma$ are corresponding
weight vectors in $\VV(\theta^\vee)$. If $\nu+\gamma\in \Delta^\vee$, then
$0\ne [v_\nu,v_\gamma]\in \g^\vee$.
\end{proof}%
Since $({\frak so}_{2p+1})^\vee={\frak sp}_{2p}$ and
$\GR{F}{4}^\vee=\GR{F}{4}$, we see that the number of
$\be$-stable commutative subalgebras is equal to
\begin{itemize}
\item  \ $2^{p-1}$ for the ${\frak sp}_{2p}$-module $\VV(\vp_2)$;
\item  \ $2$  for the ${\frak so}_{2p+1}$-module $\VV(\vp_1)$;
\item  \ $2^2$ for the $\GR{F}{4}$-module $\VV(\vp_1)$.
\end{itemize}
One should not think, however, that if $\hat\h=\h\oplus\VV$ is
a ${\Bbb Z}_2$-grading, then the number of $\be(\h)$-stable
commutative subalgebras of $\VV$ is always a power of 2.

{\it Examples.}  Straightforward computations give us the following.
\\
1. If $\hat\h={\frak gl}_{n}$ and $\h={\frak gl}_r\times {\frak gl}_{n-r}$,
then there are  $\genfrac{(}{)}{0pt}{}{n}{r}
+(n-r)\genfrac{(}{)}{0pt}{}{n-1}{r-1}$ commutative
$\be(\h)$-stable subalgebras in $\VV$.
\\[.6ex]
2. If $\hat\h={\frak sl}_{2p+1}$ and $\h={\frak so}_{2p+1}$,
then the number of such commutative subalgebras is
$2^{p+1}-1$.
\begin{rem}{Remark}   \label{g2}
In the previous exposition the case of $\g=\GR{G}{2}$ is omitted,
and the reason is that the little adjoint $\GR{G}{2}$-module is associated with
a certain ${\Bbb Z}_3$-grading.
Namely, there is an automorphism $\sigma$ of order 3 of ${\frak so}_8$
such that the fixed-point subalgebra is $\GR{G}{2}$ and
two other eigenspaces of $\sigma$ are little adjoint $\GR{G}{2}$-modules.
That is,
\[
   {\frak so}_8=\GR{G}{2}\oplus \W \oplus \W' \ ,
\]
where $\W\simeq\W' \simeq \VV(\theta_s)$ and
$[\W,\W]\subset \W'$, $[\W',\W']\subset \W$.
Still, it makes sense to speak about $\be$-stable commutative
subalgebras of $\W$. We have $\GR{G}{2}^\vee\simeq \GR{G}{2}$.
Then completely the same argument shows
that the $\be$-stable commutative subalgebras in $\W$ are in
one-to-one correspondence with the long Abelian $\be$-ideals
in $\GR{G}{2}$.
\end{rem}

\end{document}